\def \Bbb#1{{#1\kern-.800em #1}}

\def \be {\begin{equation}}
\def \eq {\end{equation}}
\def \bee {\begin{eqnarray}}
\def \eqq {\end{eqnarray}}
\def \nn {\nonumber}
\def \bea {\begin{array}{c}}
\def \eqa {\end{array}}

\newtheorem{df}{Definition}

\def \a {\alpha}
\def \b {\beta}
\def \g {\gamma}
\def \d {\delta}
\def \e {\epsilon}
\def \k {\kappa}
\def \l {\lambda}

\def \om {\omega}

\def \db {\bar{\delta}}
\def \ab {\bar{a}}
\def \bb {\bar{b}}
\def \cb {\bar{c}}
\def \etb {\bar{\eta}}
\def \deltb {\bar{\delta}}

\def \Om {\Omega}
\def \D {\nabla}
\def \Dt {\tilde{\nabla}}
\def \del {\partial}
\def \s {\;\;}
\def \G {\Gamma}
\def \Gt {\tilde{\Gamma}}
\def \Rt {\tilde{R}}
\def \P {\Phi}

\def \ib {\bar{i}}
\def \jb {\bar{j}}
\def \kb {\bar{k}}
\def \zb {\bar{z}}

\documentstyle[12pt]{article}

\begin{document}
\begin{center}
\hfill    SISSA 185/96/FM\\
\hfill    UU-HEP/96-09 \\
\hfill    q-alg/9612031\\
\vskip .1in
{\large \bf Poisson Algebra of Differential Forms}
\footnote{This work was supported in part by the Director, Office of Energy 
Research, Office of High Energy and Nuclear Physics, Division of High 
Energy Physics of the U.S. Department of Energy under Contract 
DE-AC03-76SF00098 and in part by the National Science Foundation under
grant PHY-9514797.}
\vskip .2in
Chong-Sun Chu\footnote{e-mail address: cschu@sissa.it}\\
{\em SISSA, International School for Advanced Studies\\
Via Beirut 2-4, 34013 Trieste\\
Italy\\}
\vskip .1in
Pei-Ming Ho\footnote{e-mail address: pmho@mail.physics.utah.edu}\\
{\em Physics Department\\
University of Utah\\
Salt Lake City, UT84102\\
USA\\}
\end{center}
\begin{abstract}

We give a natural definition of a 
Poisson Differential Algebra. Consistency
conditions are formulated in geometrical terms. 
It is found that one
can often locally put the Poisson structure
on the differential calculus
in a simple canonical form by a
coordinate transformation. 
This is in  analogy with  the standard Darboux's
theorem for symplectic geometry. 
For certain cases there exists a realization of the exterior
derivative through a certain canonical one-form. 
All the above 
are carried out similarly for the case of a
complex Poisson Differential Algebra. 
The case of one complex dimension
is treated in detail and interesting features are noted. 
Conclusions are made in the last section.

\end{abstract}

\newpage
\renewcommand{\thepage}{\arabic{page}}
\setcounter{page}{1}

\section{Introduction}

The motivation of the study of Poisson algebra of
differential forms is two-fold.
On the one hand we want to understand differential calculi
on quantum spaces from their Poisson limits;
on the other hand it is by itself of mathematical interest
to generalize the Poisson algebra to the differetial calculus.

In Sec.2 we give the definition of Poisson differential algebra,
which is motivated by the Poisson limit of quantum 
differential calculi \cite{CHZ,CHZ1,CHZ2,CHZ3}.
Although the definition seems to be quite natural,
the consistency conditions turn out to be very restrictive.
We will see in Sec.3 that on a local symplectic patch
the Poisson structure can always be brought to a canonical form
and the possible Poisson structures are 
labelled by a finite number of parameters, among them, for example,
are the $r$-matrices of the classical Yang-Baxter equation.
Equivalent
Poisson Differential Algebras have their
$r$-matrices related to each other
by similarity transformations.
It is shown in Sec.4 that when the torsion has a zero
there exists a one-form realization of the
exterior derivative \cite{CHZ,CHZ1,CHZ2,CHZ3}.
We extend the definition and results in Secs.2-4 to
the case of complex manifolds in Sec.5.
We notice interesting cases where there is a natural K\"{a}hler form
which gives a covariantly constant metric and which is closely
related to the one-form realization of the exterior derivative.
In Sec.6 we focus on one-dimensional complex manifolds
and find even stronger results than in the general case.
The only allowed metrics are those of a plane, a sphere
and a Lobachevskian disk.

\section{Definition of Poisson Differential  Algebras}

\begin{df} \label{realpois}
For a differetial calculus $\Om$ generated by
the coordinates $\{x^{\a}\}$ and differential forms $\{dx^{\a}\}$
on a differentiable manifold,
a bilinear map $(\cdot,\cdot):\Om\otimes\Om\rightarrow\Om$
is called a {\em Poisson structure} if for arbitrary $f,g,h\in\Om$
\begin{enumerate}
\item It has the symmetry
\be
(f,g)=(-1)^{p(f)p(g) +1}(g,f),
\eq
where the parity $p(f)=0$ if $f$ is even and $p(f)=1$ if $f$ is odd.
\item The graded Jacobi identity holds
\be
(f,(g,h))+(-1)^{p(f)(p(g)+p(h))}(g,(h,f))+(-1)^{p(h)(p(f)+p(g))}(h,(f,g))=0.
\eq
\item Elements in $\Om$ act as derivation via the map on other elements
\be
(f,gh)=(f,g)h+(-1)^{p(f)p(g)}g(f,h).
\eq
\item The Leibniz rule holds for the exterior derivative $d$
\be\label{Leibniz}
d(f,g)=(df,g)+(-1)^{p(f)}(f,dg).
\eq
\item $(f,g)$ is a form of degree equal to
the sum of the degrees of $f$ and $g$.
\end{enumerate}
The differential algebra equipped with such a Poisson structure
 is called a {\em Poisson differential algebra}.
\end{df}

The first three properties can be derived if the Poisson structure is 
obtained from the limit of a graded commutator.
They also resemble the properties of
a graded universal enveloping algebra of a Lie algebra.
The last two properties are motivated by the quantization of
the differential calculus for a quantum space
\cite{CHZ,CHZ1,CHZ2,CHZ3,Wor3,Z1}.
By studying this Poisson differential algebra we will be able to
reach a better understanding of the quantum differential calculus.
It is remarkable that this simple definition
of a Poisson structure on a differential calculus naturally leads
to geometrical notions such as connection, torsion and curvature. 

Assuming that the matrix $P^{\a\b}$ is invertible,
the Poisson algebra is completely specified by
the functions $P^{\a\b}$ and $\G^{\a}_{\b\g}$ defined by
\bee
&(x^{\a},x^{\b})=P^{\a\b}, \label{P}\\
&(x^{\a},dx^{\b})=-P^{\a\g}\G^{\b}_{\g\d}dx^{\d}. \label{Gamma}
\eqq

$P^{\a\b}$ gives the usual Poisson structure on functions.
Under the general coordinate transformation $x\rightarrow x'=x'(x)$,
it transforms as a rank-two tensor according to (\ref{P})
\be
P^{\a\b}\rightarrow P'^{\a\b}=
\frac{\del x'^{\a}}{\del x^{\g}}\frac{\del x'^{\b}}{\del x^{\d}}P^{\g\d}.
\eq
The additional data required for a Poisson structure on $\Om$
is given by $\G^{\a}_{\b\g}$,
which transforms as a connection as a result of (\ref{Gamma})
\be
\G^{\a}_{\b\g}\rightarrow \G'^{\a}_{\b\g}=
\frac{\del x^{\k}}{\del x'^{\b}}\frac{\del x^{\l}}{\del x'^{\g}}
(\frac{\del x'^{\a}}{\del x^{\d}}\G^{\d}_{\k\l}-
\frac{\del^2 x'^{\a}}{\del x^{\k} \del x^{\l}}).
\eq
We define the {\em connection one-form} for a Poisson algebra by
\be\label{G}
\G^{\a}_{\b}=dx^{\g}\G^{\a}_{\g\b}.
\eq
Note that
\be\label{Gt}
\Gt^{\a}_{\b}=\G^{\a}_{\b\g}dx^{\g}
\eq
also transforms like a connection one-form
but is different from the one above because
in general the torsion
\be \label{torsion}
T^{\a}_{\b\g}=\G^{\a}_{\b\g}-\G^{\a}_{\g\b}
\eq
does not vanish.

For a tensor $U^{\a \cdots}_{\b \cdots}$, the covariant
derivative $\D_\d U^{\a \cdots}_{\b \cdots}$ is defined by
\be 
 \D_\d U^{\a \cdots}_{\b \cdots} = \del_\d U^{\a \cdots}_{\b \cdots}
+ \G^\a_{\d \mu} U^{\mu \cdots}_{\b \cdots}
-U^{\a \cdots}_{\mu \cdots}\G^\mu_{\d \b},
\eq
and similarly for the covariant derivative $\tilde{\D}$
with respect to $\tilde{\G}$.

We will see later that the torsion plays an important role
in the Poisson algebra.

\section{Solutions of Poisson Algebras}

The Poisson structure functions $P^{\a\b}$ and $\G^{\a}_{\b\g}$
are greatly constrained by the Jacobi identities, the Leibniz rule
and the other properties of a Poisson differential algebra.
We will show in this section that on a symplectic local patch
of the manifold the Poisson structure is completely characterized by
a finite number of parameters.

The algebra of differential calculus on a local patch
is generated by the functions $\{x^{\a}\}$ and one-forms $\{dx^{\a}\}$.
We shall examine all Jacobi identities and Leibniz rules applied
to all generators.

First of all $P^{\a\b}$ has to satisfy the Jacobi identity
for three functions
\be\label{PP}
\sum_{(\a,\b,\g)}P^{\a\d}\frac{\del}{\del x^{\d}}P^{\b\g}=0,
\eq
where $(\a,\b,\g)$ indicates cyclic permutations of $\a,\b,\g$.

In the following we assume that $P^{\a\b}$ is invertible
with the inverse $P_{\a\b}$
($P_{\a\g}P^{\g\b}=P^{\b\g}P_{\g\a}=\d_{\a}^{\b}$)
giving the symplectic structure.
Eq.(\ref{PP}) can be written as
\be
\sum_{(\a,\b,\g)}\del_{\a}P_{\b\g}=0,
\eq
where $\del_{\a}$ denotes $\frac{\del}{\del x^{\a}}$,
or equivalently,
\be
dP=0,
\eq
where $P$ is the symplectic two-form
\be \label{sympl}
P=\frac{1}{2}P_{\a\b}dx^{\a}dx^{\b}.
\eq

Using Eq.(\ref{PP}), one can show that
the Jacobi identity for two functions and one one-form
requires that the Riemann curvature vanishes
\be\label{R=0}
R^{\a}_{\s\b}=\frac{1}{2}R^{\a}_{\s\b\g\d}dx^{\g}dx^{\d}=0,
\eq
where
\be
R^{\a}_{\s\b}=d\G^{\a}_{\b}+\G^{\a}_{\g}\G^{\g}_{\b}.
\eq
The vanishing of the curvature implies 
that the connection is a pure gauge.
Hence there exists a square matrix of functions $M^{\a A}$
such that
\be\label{GM}
\G^{\a}_{\b\g}=M^{\a A}\del_{\b}M_{A\g},
\eq
where $M_{A\a}$ is the inverse of $M^{\a A}$.
Note that (\ref{GM}) defines $M^{\a A}$
up to an invertible constant matrix $N^A_B:
M^{\a A}$ and $M^{\a B} N^A_B$ give the same connection.

The Leibniz rule for two functions
\be
d(x^{\a},x^{\b})=(dx^{\a},x^{\b})+(x^{\a},dx^{\b})
\eq
implies that
\be \label{Leib1}
\Dt P^{\a\b}=dP^{\a\b}+P^{\a\g}\Gt^{\b}_{\g}+\Gt^{\a}_{\g}P^{\g\b}=0.
\eq
Plugging (\ref{GM}) into this equation we find
\be
\del_{\a}G^A_{\b}-\del_{\b}G^A_{\a}=0,
\eq
where $G^A_{\a}=P_{\a\b}M^{\b A}$.
This implies that locally we can find functions $\{\P^A\}$
such that
\be
G^A_{\a}=\del_{\a}\P^A.
\eq
Hence
\be \label{MPdP}
M^{\a A}=P^{\a\b}\del_{\b}\P^A.
\eq

The Leibniz rule for one function and one one-form
follows from (\ref{Leib1}).

Using Eqs.(\ref{PP}) and (\ref{R=0}),
one can show after considerable calculations that
the Jacobi identity for one function and two one-forms
implies
\be \label{DPR}
\D_{\sigma}(P^{\a\g}\Rt^{\b}_{\s\g\k\l})=0.
\eq

With the Leibniz rule (\ref{Leibniz})
the Jacobi identity for three one-forms follows.
It is easy to convince oneself that Eqs.(\ref{PP}), (\ref{R=0}), (\ref{Leib1})
and (\ref{DPR}) are all the conditions on $P^{\a \b}$ and $\G^{\a}_{\b \g}$
that a Poisson differential algebra has to satisfy.

Now, we are going to prove that for any Poisson differential algebra
with a symplectic two-form (\ref{sympl}),
there is a preferred local coordinate chart 
in which the Poisson structure is a quadratic form
of the local coordinates (Eq.(\ref{PAB}) below) 
and the connection coefficients can be expressed in 
terms of the same Poisson structure in a simple manner (Eq.(\ref{GP}) below).
This is in some sense an analogue of the Darboux' theorem for 
symplectic geometry \cite{Wein1}.

Since $M$ and $P$ are both invertible matrices,
$G$ is also invertible.
Hence by Eq.(\ref{MPdP}) $\{\P^A\}$ can be used as local coordinates.
In this new coordinate system $\{\P^A\}$,
we have
\bee
&M^{AB}=P^{AB}, \\
&\G^A_{BC}=P^{AD}\del_B P_{DC}, \label{GP} \\
&T^A_{BC}=P^{AD}\del_D P_{BC},
\eqq
by Eqs.(\ref{torsion}),(\ref{GM}) and (\ref{MPdP}).

It is interesting that in this coordinate system
everything can be expressed in terms of
$P^{AB}=(\P^A,\P^B)$ alone.
Note that due to the ambuigity in the matrix $M^{\a A}$, one can perform
the  transformation 
\be \label{phi_transf}
\P^A \rightarrow N^A_B \P^B +V^A
\eq
on the coordinate $\P^A$ 
without changing the connection coefficient $\G^\a_{\b \g}$.
Here $N^A_B$ and $V^A$ are constants.

Define a new basis of one-forms
\be \label{e_A}
e_A \equiv M_{A\a}dx^{\a}=P_{AB}d\P^B.
\eq
The connection vanishes in this basis
\be \label{nabla-del}
\D^A=P^{AB}\del_B
\eq
and the torsion is
\be\label{tor}
T_A^{\s BC}=\del_A P^{BC}.
\eq
Since we use the same kind of indices ($A,B,C,\cdots$)
for both bases $\{e_A\}$ and $\{d\P^A\}$,
one has to keep in mind which basis one is dealing with
for a tensor with indices ($A,B,C,\cdots$). We will say it explicitly
whenever there could be a possible confusion.

It can be checked that
\be \label{ee}
(e_A,e_B)=-\Rt_{AB}=-\frac{1}{2}\Rt_{AB}^{CD}e_C e_D,
\eq
where $\Rt^A_{\s BCD}$ are the compoments of the curvature tensor in the 
$\{d\Phi^A\}$ bases,
\be \Rt^A_{\s BCD}= G^A_{\a} G^{\b}_B G^{\g}_C G^{\d}_D \Rt^{\a}_{\s\b\g\d},
\eq
and
$\Rt_{AB}^{CD}$ comes from $\Rt^A_{\s BCD}$ by raising and lowering indices
with the $P^{AB}$,
\be
\Rt_{AB}^{CD}=P_{AE}P^{CF}P^{DG}\Rt^E_{\s BFG}.
\eq
It satisfies
\be
\Rt_{AB}^{CD}=\Rt_{BA}^{CD}=-\Rt_{AB}^{DC}.
\eq

The basis $\{e_A\}$ worths special treatment because
the Poisson brackets between
$e_A$ and functions vanish
\be
(e_A,x^{\a})=0.
\eq
Calculations in this basis could be
much simpler than those in others.
For example, to check the Jacobi identity for
one function and two one-forms,
it is sufficient to check
\be
(x^{\a},(e_A,e_B))=0,
\eq
which immediately gives
\be
\Rt_{AB}^{CD}=\mbox{constants}.
\eq
The same result can be obtained by expressing Eq.(\ref{DPR})
in the basis $\{e_A\}$ using (\ref{nabla-del}).

Using (\ref{R=0}) and the following two identities:
\bee
&\sum_{(\b,\g,\d)}(R^{\a}_{\s\b\g\d}+\D_{\b}T^{\a}_{\g\d}
-T^{\a}_{\b\k}T^{\k}_{\g\d})=0, \\
&\Rt^{\a}_{\s\b\g\d}-R^{\a}_{\s\b\g\d}=-\D_{\g}T^{\a}_{\d\b}
-\D_{\d}T^{\a}_{\b\g}+T^{\a}_{\b\k}T^{\k}_{\g\d}
+T^{\a}_{\g\k}T^{\k}_{\d\b}+T^{\a}_{\d\k}T^{\k}_{\b\g},
\eqq
we find
\be
\Rt^{\a}_{\s\b\g\d}=\D_{\b}T^{\a}_{\g\d}.
\eq
In the basis $\{e_A\}$ it is
\be
\Rt_{AB}^{CD}=\del_B T_A^{\s CD},
\eq
which can be easily solved
\be
T_A^{\s CD}=\Rt_{AB}^{CD}\P^B+f_A^{CD},
\eq
where $f_A^{CD}=-f_A^{DC}$ are constants.
Now $P^{AB}$ can be solved from (\ref{tor})
\be \label{PAB}
P^{AB}=(\P^A,\P^B)
=\frac{1}{2}\Rt^{AB}_{CD}\P^C\P^D+f^{AB}_C\P^C+g^{AB},
\eq
where $g^{AB}=-g^{BA}$ are constants.

We will call $\{\P^A\}$ the canonical coordinate system.
If the torsion vanishes, the canonical coordinates
will coincide with the Darboux coordinates
up to the transformation (\ref{phi_transf}),
and the Poisson structure between function and forms
will be trivial.

It is remarkable that
in the canonical coordinate system
all information about
a Poisson algebra with invertible $P^{A B}$
is encoded in the constants
$\Rt^{AB}_{CD}, f^{AB}_C$ and $g^{AB}$.

It is interesting that the curvature $\Rt^{A B}_{C D}$ must
satisfy the classical Yang-Baxter equation. This can be
seen as a result of the consistence of the Jacobi
identity.
Using the tensor product notation
\bee
&\Rt_{12}= \Rt \otimes 1, \\
&(\Rt_{12})^{AB}_{CD}= \Rt^{AB}_{CD},
\eqq
(\ref{PAB}) can be written as 
\be 
(\Phi_1, \Phi_2)= \frac {1}{2}\Rt_{12} \P_1\P_2+f_{123} \P_3+g_{12},
\eq
where matrix multiplication is implied.
The Jacobi identity of three functions implies
\bee 
0&=&\sum_{(1,2,3)}(\Phi_1,(\Phi_2,\Phi_3)) \nn \\
&=& - \frac{1}{4} ([\Rt_{12},\Rt_{13}]+[\Rt_{12},\Rt_{23}]
+[\Rt_{13},\Rt_{23}] ) \Phi_1 \Phi_2 \Phi_3  \nn \\
& & + \frac{1}{2}\sum_{(1,2,3)}
(\Rt_{23}(f_{124}\P_3\P_4+f_{134}\P_2\P_4) + f_{234}\Rt_{14}\P_1\P_4)\nn\\
& & + \sum_{(1,2,3)}
(\frac{1}{2}(\Rt_{13}+\Rt_{23})g_{12} \Phi_3+ f_{234} f_{145} \P_5)\nn \\
& & + \sum_{(1,2,3)}(f_{234} g_{14}),
\eqq
where $(1,2,3)$ stands for cyclic permutation of 1,2 and 3. 
Vanishing of the $\Phi_1 \Phi_2 \Phi_3$ term gives 
the classical Yang-Baxter equation,
\be \label{CYBE}
[\Rt_{12},\Rt_{13}]+[\Rt_{12},\Rt_{23}]
+[\Rt_{13},\Rt_{23}] =0.
\eq
Other conditions can be obtained from the vanishing
of the constant, $\P$ and $\P \P$ terms.
Explicitly they are
\bee
&\sum_{(A,B,C)}(2\Rt^{AB}_{FD}f^{CF}_E+f^{AB}_{F}\Rt^{CF}_{DE})=0, \\
&\sum_{(A,B,C)}(\Rt^{AB}_{ED}g^{CE}+f^{AB}_{E}f^{CE}_{D})=0,\label{ff} \\
&\sum_{(A,B,C)}f^{AB}_{D}g^{CD}=0.
\eqq
We will  see in the next section that solutions with $f^{AB}_C$=0
is of particular interest.

\section{Properties of Poisson Algebras}

Define
\be \label{xi}
\xi=-e_A\P^A.
\eq
Then
\be
(\xi,f)=df
\eq
for an arbitrary function $f$.
One can check that 
\be (\xi, dx^\a)= -\frac{1}{2} M^{\a A} f^{C D}_A e_C e_D. \eq
Therefore
the one-form $\xi$ will become
a one-form realization of the exterior derivative
\cite{CHZ,CHZ1,CHZ2,CHZ3}
\be
(\xi,\om)=d\om
\eq
for an arbitrary differential form $\om$ if
$f^{AB}_C=0$.

Note that sometimes it is possible to make
a transformation (\ref{phi_transf})
so that in the new coordinate system $f^{AB}_{C}=0$
and a one-form realization (\ref{xi}) exists.
The existence of such a $\P$-coordinate system
is equivalent to the existence of a zero
of the torsion tensor.
This can be seen from the following expression of $P^{AB}$:
\be \label{RTP}
P^{AB}=\frac{1}{2}\Rt^{AB}_{CD}\P^C\P^D+T^{AB}_{C}(0)\P^C+P^{AB}(0),
\eq
where $T^{AB}_{C}(0)$ and $P^{AB}(0)$ are the torsion tensor
and Poisson tensor evaluated at the origin $\P^A=0$.
Since the transformation (\ref{phi_transf}) can translate
the origin to any other point on the patch,
if the torsion has a zero
one can always choose the new origin
to reside on that point.
In addition, Eq.(\ref{RTP}) shows that the maps from the space
of all $\P$-coordinate systems,
where an element is labelled by $(N^A_B,V^A)$,
to the coefficients $\Rt^{AB}_{CD}, f^{AB}_{C}, g^{AB}$
of the Poisson tensor
are the same as the composite of a rotation by $N^A_B$ with
the tensors $\Rt^{AB}_{CD}, T^{AB}_{C}, P^{AB}$
as maps from the local patch to the values of those tensors.

While in general
\be
d\xi=(\frac{1}{2}f^{AB}_C\P^C+g^{AB})e_A e_B,
\eq
when $f^{AB}_{C}=0$
it is tempting to interpret $d\xi$ as the length element squared,
$g^{AB}$ as the metric and $e_A$ as the vielbein.
The tensor $g^{A B}$ defined this way is covariantly constant
with respect to the connection $\G$.
But since $g^{AB}$ is antisymmetric and $d\xi$, $e_A$ are forms,
in fact $d\xi$ is more like the K\"{a}hler two form
for a complex manifold.

As a result of (\ref{nabla-del}),
a metric $h_{\a \b}$ is covariantly constant,
\be \label{nabla-h}
\nabla h_{\a\b}=0
\eq
if and only if its components $h^{A B}$ in the $\{e_A\}$ basis are constant.

Before going into the next topic,
we remark that if a covariantly constant metric
$h$ is desired on the manifold,
then (\ref{nabla-h}) and (\ref{Leib1}) determine the connection to be
\be
\G^{\a}_{\b\g}=\frac{1}{2}P_{\b\d}h_{\g\e}
(h^{\e\k}\partial_{\k}P^{\a\d}+h^{\a\k}\partial_{\k}P^{\d\e}
-h^{\d\k}\partial_{\k}P^{\e\a}+P^{\e\k}\partial_{\k}h^{\a\d}
-P^{\a\k}\partial_{\k}h^{\d\e}-P^{\d\k}\partial_{\k}h^{\e\a})
\eq
without using any other relations.
This result may be relevant to certain generalizations
of Riemannian geometry as
Eq.s (\ref{nabla-h}) and (\ref{Leib1}) relate $h$ and $P$ to
a covariantly constant complex Hermitian tensor.

\section{Complex Poisson Differential Algebras}

All of the above can be easily specialized for the case of a complex
differentiable manifold,
\begin{df}
For a differetial calculus $\Om$ generated by
the coordinates $\{ x^i, x^{\ib} \}$ and differential forms 
$\{dx^{i}, dx^{\ib} \}$ on an $N$ dimensional complex manifold,
the map $(\cdot,\cdot):\Om\otimes\Om\rightarrow\Om$
is called a {\em Poisson structure} on $\Om$ if
for arbitrary $f,g,h\in\Om$, it
satisfies properties 1 to 3 of Definition \ref{realpois} and
also

4. The Leibniz rule holds for the holomorphic and 
anti-holomorphic exterior derivatives $\d, \db$:
\bee \label{Leib-cpx}
&&\d (f,g)=(\d f,g)+(-1)^{p(f)}(f,\d g),\\
&&\db (f,g)=(\db f,g)+(-1)^{p(f)}(f,\db g).
\eqq

5. The holomorphic and anti-holomorphic degrees of the form $(f,g)$ are 
equal to the sums of those of 
$f$ and $g$.

6. Hermiticity: 
\be (f,g)^*= (-1)^{p(f)p(g)}(g^*,f^*),
\eq
where the $*-$operation on $\Om$ is the complex conjugation
together with an inversion of ordering in a product:
\be (f g)^* = g^* f^* = (-1)^{p(f)p(g)} f^* g^*. \eq

A  complex differential algebra equipped with such a Poisson 
structure is called  a 
{\em complex Poisson differential algebra}.
\end{df}

Using the notation $\a \in I \cup \bar{I}$,
where $I =\{i,j,k,\cdots \}$ is the
set of holomorphic indices and $\bar{I} =\{ \ib, \jb, \kb, \cdots \}$ is
the set of anti-holomorphic indices. 
Assuming that the matrix $P^{\a \b}$ is invertible, the 
complex Poisson differential algebra is defined by the
functions $P^{\a \b}, \G^\a_{\b \g}$ as in (\ref{P}) and (\ref{Gamma}),
where $\{x^\a\} =\{ x^i, x^{\ib} \}.$

Vanishing of the curvature $R^\a_{\s\b\g\d} =0$ implies that
connection is of the form (\ref{GM}),
which can be solved as
\be \label{MG}
M_{A\a}=M_{A\b}(0)(Pexp\int_0 \G)^{\b}_{\a},
\eq
where $M_{A\a}(0)$ is the value of $M_{A\a}$
at a certain point where
the path for the line integral $\int_0$ begins,
$P$ stands for path ordering
and $\G$ is the matrix of one-forms $\G^{\a}_{\b}$.
$M_{A\a}(0)$ can be any invertible constant matrix.
We choose it to be block-diagonal
\be \label{blockM}
M_{A\a}(0)=\left(
           \begin{array}{cc}
             M_{ai} & 0 \\
             0 & M_{\ab\ib}
           \end{array}
           \right),
\eq
where we split the index $A$ into $a$ and $\ab$,
with no implication on the holomorphicity.

The new input property 5 implies that the following 
components of the connection vanishes
\be \label{G=0} 
\G^j_{\a {\kb}} = \G^{\jb}_{\a k} =0
\eq
so that $\G$ is also block-diagonal.
It follows from (\ref{MG}) that $M_{A\a}$ must be block-diagonal as well.

The basis of one-forms (\ref{e_A}) now splits
into the holomorphic part $\{e_a\}$
and the antiholomorphic part $\{e_{\ab}\}$, where
\be 
e_a = M_{ai} dx^i, \quad 
e_{\ab} = M_{\ab\ib} dx^{\ib},
\eq
satisfy
\be 
(e_a, x^i)=(e_a, x^{\ib}) =0.
\eq

For arbitrary nondegenerate Hermitian constant matrix $h^{a\bb}$,
$h^{a\bb}e_{a}e_{\bb}$ is an admissible K\"{a}hler
two-form \cite{nakahara} which gives a covariantly constant metric and
is central to all functions.

It follows from property 5 and (\ref{ee}) that
$\Rt^{AB}_{CD}$ vanishes unless
the number of barred (unbarred) superscripts
equals the number of barred (unbarred) subscripts.
For example,
\be \Rt^{d \cb}_{a b} = \Rt^{d \cb}_{\ab \bb} =0. \eq

Given the matrix $M^{\a A}$
one can find $\P^A$ up to additive constants from (\ref{MPdP}).
$\{\P^A\}$ is split into $\{\phi^a\}$ and $\{\phi^{\ab}\}$.
In terms of these functions $\{ \phi^a, \phi^{\ab} \}$,
the Poisson structure is particularly simple.
For example,
\be \label{Pab}
P^{ab} = (\phi^a, \phi^b) =
\frac{1}{2}\Rt^{ab}_{cd}\phi^c \phi^d +f^{ab}_c \phi^c + g^{ab}
\eq
and
\be \label{Pabb}
P^{a \bb} = (\phi^a, \phi^{\bb}) =
\Rt^{a \bb}_{c \bar{d}}\phi^c \phi^{\bar{d}} +f^{a \bb}_c \phi^c +
f^{a \bb}_{\cb} \phi^{\cb}+ g^{a \bb}.
\eq

In general $(\phi^{a})^*$ is different from $\phi^{\ab}$.
But using the Hermiticity of the Poisson structure (property 6)
and (\ref{MPdP}) one sees that by choosing $M_{A\a}(0)$ to be
of the form
\be \label{M*}
M_{A\a}(0)=\left(
           \begin{array}{cc}
              m_{ai} & 0 \\
              0 & m_{\ab\ib}
           \end{array}
           \right)
\eq
where $m_{\ab\ib}$ is defined as
\be m_{\ab\ib} \equiv -m_{a i}^* \eq
and $m$ is some invertible matrix,  we get
\be
(\del_{i}\P^a)^*=\del_{\ib}\P^{\ab},
\quad (\del_{i}\P^{\ab})^*=\del_{\ib}\P^a
\eq
so that $\P^A$ can be chosen to satisfy
\be \label{phi-phi*}
\phi^{a*}=\phi^{\ab}.
\eq
In this case,
\be \label{Rt*}
(\Rt^{AB}_{CD})^*=-\Rt^{\bar{A}\bar{B}}_{\bar{C}\bar{D}},
\eq
where $\bar{A}=(\ab,\bar{\ab})=(\ab,a)$, and
\be
(e_a)^*= -e_{\ab}.
\eq

The coordinates $(\phi^a,\phi^{\ab})$
are fixed up to transformations
\be \label{ntransf}
\phi^a\rightarrow n^a_b\phi^b+v^a
\eq
and its complex conjugation. Here $n^a_b$ is an  arbitary invertible
constant matrix.

We have shown that given a complex Poisson differential algebra
one can always find a coordinate system satisfying
all the properties
(\ref{Pab}), (\ref{Pabb}), (\ref{phi-phi*}) and (\ref{Rt*}).

From the point of view of classifying all complex Poisson
differential algebras,
one is interested in how many different complex structures
there can be for a certain Poisson structure given
in terms of $\phi^a, \phi^{\ab}$.
On this matter we only comment
that the condition (\ref{blockM}) restricts the transformation
$x^{\a} \rightarrow \Phi^A$ to satisfy the partial differential equation
\be
\frac{\partial x^{\ib}}{\partial \phi^a} P^{a b} +
     \frac{\partial x^{\ib}}{\partial \phi^{\ab}} P^{{\ab} b} =0
\eq
and inequivalent complex structures correspond to
inequivalent holomorphic classes of solutions to this equation.
The coordinates $(x^i,x^{\ib}=x^{i*})$ and
$(y^i,y^{\ib}=y^{i*})$ are said to be in the same holomorphic class if
${\partial x^i}/{\partial y^{\ib}}=0$.

In the coordinate system $(\phi^a,\phi^{\ab}=\phi^{a*})$,
the one-forms 
\be \eta \equiv -e_a \phi^a, \quad \etb \equiv - e_{\ab} \phi^{\ab}
\eq
are holomorphic and antiholomorphic respectively.
They satisfy
\be
\eta^* = -\etb \label{e0} 
\eq
and
\bee 
(\eta,x^i)=dx^i &,&(\eta,x^{\ib})=0, \label{e1} \\
(\etb,x^{\ib})=dx^{\ib} &,&(\etb,x^i)=0. \label{e2}
\eqq

For those solutions with  $f^{A B}_C=0$
(possibly after a transformation (\ref{ntransf})),
we have the one-form realization for
$\delta$ and $\deltb$:
\be \delta \om =(\eta, \om), \quad \deltb \om =(\etb, \om),
\eq
for any $\om \in \Om$.
For this case, the 2-form
\be \label{K}
K \equiv \db \eta = g^{a \bb} e_a e_{\bb} = \d \etb
\eq
is central in the Poisson bracket
\be \label{Kcen} (K, x^\a)=(K, dx^\a)=0 \eq
as a consequence of (\ref{e1}) and (\ref{e2}).
It is not hard to check that
\be
\d \eta = g^{a b} e_a e_b, \quad
\db \etb = g^{\ab \bb} e_{\ab} e_{\bb}.
\eq

If $g^{a b}= g^{\ab \bb}=0$
\be
\d \eta = \db \etb =0,
\eq 
then at least locally 
\be
\etb =\db V
\eq 
for certain function $V$, and
\be K=\d \db V
\eq
shows that it is a K\"{a}hler manifold.

We will see in the next section that in the case of one dimensional
complex Poisson differential algebra,
the admissible K\"{a}hler metrics are those 
of constant curvature spaces.

\section{One Dimensional Case}

For the case of one complex dimension, the Poisson structure is completely
determined by the two functions $P, S$ which appear in
\be (z,\zb) = P, \quad (z,dz) =Sdz,
\eq
where $P$ is real and $S$ is holomorphic
due to properties 4 and 6.
All other relations of the Poisson structure 
can be derived by using the Leibniz rule.
The Jacobi identity
\bee 0 &=& ((z,\zb),dz) + ((\zb, dz), z)+ ((dz, z), \zb) \\ \nn
 &=& -P^2 [\del(P^{-1} S) - \del (P^{-1} \bar{\del} P)]
\eqq
implies that 
\be S= (\bar{\del} +\bar{\a}) P \eq
for some $\bar{\a} =\bar{\a}(\zb)$.
Now, under a holomorphic coordinate transformation 
$z \rightarrow z'=f(z)$,
\bee
&P' =({\del f}) (\overline{ {\del f}}) P, \\
&S' = ({\del f}) S, \\
&\bar{\del'} = ( \overline{ {\del f}})^{-1} \bar{\del}
\eqq
and so
\be
\a'= ({\del f})^2 [(\del f) \a - {{\del}^2 f}].
\eq
In particular, one can pick a coordinate transformation to make 
$\a'=0$. Omitting the primes,
we have just  shown that one can always make a holomorphic
coordinate transformation and arrive at the relation 
\be S=\bar{\partial} P,
\eq
in a particular coordinate system.
Because $P$ should be real, it must be
\be P=a z \zb +b z +\bb \zb +c,
\eq
where the constant coefficients $a,c$ are real and $b$ is complex.

In the case of higher dimensions, the coordinate system
$\{\phi^a, \phi^{\ab}\}$, which gives
the quadratic form of the Poisson structure,
are generically not holomorphic nor antiholomorphic
functions of $z^i, \zb^{\ib}$ in which the complex structure
of the manifold is defined.
It is remarkable that in the case of one dimension,
holomorphic transformation
is enough to bring oneslf to the canonical coordinate system.

For the metric to be covariantly constant,
the only admissible K\"{a}hler form (up to normalization) is
\be K=dz d\zb h  \eq
with the metric $h$ given by
\be h=P^{-2}. \eq
$K$ is central to functions and forms.

Performing a fractional transformation
\be z \rightarrow \frac{\a z'+\b}{\g z' +\d},
\eq
we have
\be 
K= dz' d\zb ' h^{'},
\eq
with
\be h'= \frac {{\mid \a \d -\b \g \mid}^2}
{(A z' \zb' +B z' +\bar{B} \zb'+C)^2}
\eq
and 
\be
\left(\begin{array}{cc}
         A       & B \\
         \bar{B} & C
      \end{array}\right)=
\left(\begin{array}{cc}
         \a & \g \\
         \b & \d
      \end{array}\right)
\left(\begin{array}{cc}
         a   & b \\
         \bb & c
      \end{array}\right)
\left(\begin{array}{cc}
         \bar{\a} & \bar{\b} \\
         \bar{\g} & \bar{\d}
      \end{array}\right).
\eq

One can always choose the fractional transformation with 
$\a \d -\b \g \neq 0$ such that $B=0$ and so
\be h'= \frac {{\mid \a \d -\b \g \mid}^2}{(A z' \zb' +C)^2}.
\eq 
This is nothing but the most general form of a metric
for a one dimensional Hermitian surface 
with constant Gaussian curvature (for the Levi-Civita connection).
For $A=0$ or $C =0$ (do $z''=1/z'$ in this case), 
$h'$ gives the metric of a plane.
For $A \neq 0$, the metric is that of a sphere (resp. Lobachevskian plane)
if $A$ and $C$ are of the same sign (resp. opposite sign). 
All three cases are actually K\"{a}hler manifolds.
It is interesting to note that all connected, simply connected
one-dimensional complex manifolds are biholomorphic to
one of the three cases or a quotient of them over an isometric automorphism.

If $a\neq 0$ one can always make a translation for $z,\zb$
so that $b=\bar{b}=0$ in the new coordinate system.
It can be checked that in the new coordinate system
the holomorphic and antiholomorphic exterior derivatives
are realized by
\be
\eta=-\zb P^{-1}dz, \quad \bar{\eta}=z P^{-1}d\zb,
\eq
which gives an admissible K\"{a}hler form described above by (\ref{K}).

\section{Conclusion}

We found in this paper that the generalization of the Poisson structure
to the differential calculus,
or the Poisson limit of the differential calculus
on a quantum space, has interesting geometrical structure in itself.
We proved an analogue of Darboux's theorem
for the present case of a Poisson
Differential Algebra.
We also studied the compatibility of the complex structure with
the Poisson structure on complex manifolds.
In the case of one complex dimension the natural metrics
are the ones whose Levi-Civita connections give constant curvatures.
It is natural to ask if one can generalize the results
for one complex dimension to higher dimensions.
It will also be interesting to see if there exist a canonical
procedure to quantize a Poisson differential algebra.
We leave these questions for a future publication.


\section{Acknowledgements}
We are very grateful to Bruno Zumino
for his essential help in the form of 
discussions, advice
and computations.  
This work would not have been possible without his help.

\end{document}